\documentclass{amsart}

\usepackage{enumitem}
\usepackage{amsfonts, amssymb, amsthm, amsmath, calc, cancel, cite,color, eucal, graphics, graphicx, hyperref, latexsym, mathdots, multirow, pgfplots, theoremref, tikz,tikz-cd, url}
\usepackage{accents}
\usepackage{caption} 
\numberwithin{equation}{section}

\theoremstyle{plain}
\newtheorem{theorem}{Theorem}[section]
\newtheorem{lemma}[theorem]{Lemma}
\newtheorem{proposition}[theorem]{Proposition}
\newtheorem{corollary}[theorem]{Corollary}

\theoremstyle{definition}
\newtheorem{definition}[theorem]{Definition}

\newtheorem{remark}[theorem]{Remark}
\newtheorem{example}[theorem]{Example}

\renewcommand{\phi}{\varphi}
\renewcommand{\emptyset}{\varnothing}
\newcommand{\meet}{\wedge}
\newcommand{\join}{\vee}

\newcommand{\mb}[1]{\mathbf{#1}}
\newcommand{\bb}[1]{\mathbb{#1}}

\newcommand{\mt}[1]{\mathtt{#1}}
\newcommand{\mc}[1]{\mathcal{#1}}

\newcommand{\A}{\mb{A}}
\newcommand{\B}{\mb{B}}
\newcommand{\C}{\mb{C}}
\newcommand{\D}{\mb{D}}
\newcommand{\I}{\mb{I}}
\newcommand{\M}{\mb{M}}
\newcommand{\U}{\mb{U}}

\newcommand{\comm}{\mathrm{comm}}
\newcommand{\nil}{\mathrm{nil}}

\newcommand{\bdot}{\boldsymbol{\cdot}}

\DeclareMathOperator{\Con}{Con}
\DeclareMathOperator{\Id}{Id}
\DeclareMathOperator{\glb}{glb}
\DeclareMathOperator{\DI}{DI}
\DeclareMathOperator{\Ht}{ht}

\title[Nilpotent BCK-algebras]{Nilpotent BCK-algebras}

\author[C. M. Evans]{C. Matthew Evans}
\address{Mathematics Department\\
Washington \& Jefferson College\\Washington, PA\\USA}
\urladdr{https://sites.google.com/view/mattevans}
\email{mevans@washjeff.edu}



\begin{document}

\begin{abstract}
We use the notion of pseudocommutators to define the \emph{derived ideal} of a BCK-algebra. Using the derived ideal we define the \emph{commutativization} of a BCK-algebra; this construction is functorial and the commutativization functor is left adjoint to the inclusion functor from the category of commutative BCK-algebras to the category of BCK-algebras. This yields a new proof that commutative BCK-algebras are a reflective subcategory of BCK-algebras. After this, we introduce central series and define a notion of \emph{nilpotence} for BCK-algebras and prove some properties of nilpotence. In particular, for any variety of BCK-algebras, the subclass of nilpotent algebras is a subpseudovariety, though in general not a subvariety. We also show that the class of BCK-algebras of nilpotence class at most $c$ is a subquasivariety of all BCK-algebras, and is a variety if and only if $c=1$. We close by showing that every BCK-algebra of finite height is nilpotent.
\end{abstract}

\maketitle

\section{Introduction}

As is commonly known, Boolean algebras are an algebraic formulation of classical propositional logic; in the more precise language of Blok and Pigozzi \cite{BP89}, Boolean algebras are the \textit{equivalent algebraic semantics} for classical propositional logic. There are many algebraizable logics, too many to give a full account here. In this paper, we are concerned with BCK-algebras, which generalize Boolean algebras. 

In a series of twenty-five papers appearing in the \textit{Proceedings of the Japan Academy} in 1965 and 1966, authors Arai, Imai, Is\'{e}ki, and Tanaka considered axiom systems for a variety of propositional calculi. Among the many propositional calculi considered in these papers, one is the BCK-system determined by \textit{modus ponens} together with the axioms
\begin{enumerate}
\item[(B)] $(q\rightarrow r)\rightarrow\bigl((p\rightarrow q)\rightarrow (p\rightarrow r)\bigr)$
\item[(C)] $\bigl(p\rightarrow(q\rightarrow r)\bigr)\rightarrow\bigl(q\rightarrow(p\rightarrow r)\bigr)$
\item[(K)] $p\rightarrow(q\rightarrow p)$\,.
\end{enumerate} 
Thus, BCK is a purely implicational system. This system was isolated as early as 1934 by Tarski in \cite{Tar34} (see also \cite[p.397]{Tar56}). The corresponding algebraic formulation of this system was given by Is\'{e}ki in \cite{iseki66} and Is\'{e}ki and Imai in \cite{II66}. These algebras are named in reference to the earlier combinatory logic of Sch\"{o}finkel \cite{schon24} and Curry \cite{curry30}, specifically the combinators labeled B, C, and K. 

Many of the commutative algebras of logic -- Boolean algebras, MV-algebras, Heyting algebras, Hilbert algebras, etc -- can be expressed as special cases of BCK-algebras; Iorgulescu's text \cite{iorg08(2)} provides a very detailed study of these connections. We highlight a small handful of connections to other algebraic structures. The class of bounded commutative BCK-algebras is categorically equivalent to the class of MV-algebras \cite{mundici86}, which are the algebraic semantics for many-valued logic. 
Bounded \textit{implicative} BCK-algebras are (term reducts of) Boolean algebras \cite{it78}. In 1980, Cornish demonstrated a categorical equivalence between lattice-ordered Abelian groups and so-called \textit{conical} BCK-algebras \cite{cornish80}; these are BCK-algebras satisfying a property called \textit{condition (S)} together with two additional identities. Earlier, Is\'{e}ki had shown that BCK-algebras with condition (S) are commutative semigroups with respect to a derived operation \cite{iseki79}. Lastly, we mention that certain classes of rings admit an implicative BCK-structure \cite{cornish81}. 

One key difference between BCK-algebras and the other algebras of logic mentioned above is the lack of {\it commutativity}. Specifically, the commutative law ``$x\meet y= y\meet x$'' is satisfied in any classical setting, but also some non-classical settings, including MV-algebras and Heyting algebras. However, this kind of commutativity, to be defined more precisely in Section \ref{prelims}, is not true in BCK-algebras in general. One may then wonder: given a BCK-algebra, ``how commutative'' is it? 

One way to measure this defect in finite BCK-algebras is to look at the quotient
\[\frac{\text{(\# pairs of commuting elements)}}{\text{(total \# of pairs of elements)}}\,,\] called the {\it commuting degree}, which is just the probability that two randomly selected elements commute. This probabilistic approach was used to great effect for Heyting algebras -- which are the algebraic semantics for intuitionistic logic -- in the paper \cite{BK25}. Inspired by that paper, the present author's preprint \cite{evans25} studies commuting degrees for finite BCK-algebras. Of course, the commutative law is only way to proceed; there are other equations satisfied in the classical setting that are not generally true for BCK-algebras, and one may study similarly defined probabilities for those equations. Such a study was taken up in the author's earlier paper \cite{evans23}.

The downfall to this approach is that commuting degree, as defined above, can only be computed for finite structures. In this paper, rather than using probability, we define a notion of {\it nilpotence} for BCK-algebras which gives us a measure of commutativity for the whole class of BCK-algebras. 

The paper is structured as follows: in the next section, we define BCK-algebras and discuss the elementary properties that will be needed. In Section \ref{commutators}, we introduce pseudocommutators and the derived ideal. We go on to show that $\mt{cBCK}$, the category of commutative BCK-algebras, is a reflective subcategory of $\mt{BCK}$ by explicitly constructing the reflector. In Section \ref{nilpotence}, we define central series for BCK-algebras and use them to define nilpotence. We show that nilpotence can be characterized by the behavior of both an upper central series and a lower central series. Some examples are provided. Finally, in Section \ref{nilp props}, we discuss the algebraic properties of nilpotent BCK-algebras. In particular, we show that, for any variety $\mt{V}$ of BCK-algebras, the subclass $\mt{V}_{\nil}$ of nilpotent members is a pseudovariety, though not generally a variety. We also show that $\mt{BCK}_c$, the class of BCK-algebras of nilpotence class at most $c$, is a subquasivariety of $\mt{BCK}$, and that $\mt{BCK}_c$ is a variety if and only if $c=1$. We end the paper by showing that all BCK-algebras of finite height are nilpotent.

\section{Preliminaries}\label{prelims}

\begin{definition}\label{BCK}A \emph{BCK-algebra} is an algebra $\A=\langle A; \bdot, 0\rangle$ of type $(2,0)$ such that 
\begin{enumerate}[leftmargin=1.5cm]
\item[(BCK1)] $\bigl((x\bdot y)\bdot(x\bdot z)\bigr)\bdot(z\bdot y)=0$
\item[(BCK2)] $\bigl(x\bdot (x\bdot y)\bigr)\bdot y=0$
\item[(BCK3)] $x\bdot x=0$
\item[(BCK4)] $0\bdot x=0$
\item[(BCK5)] $x\bdot y=0$ and $y\bdot x=0$ imply $x=y$.
\end{enumerate} for all $x,y,z\in A$.
\end{definition}

We give here two motivating examples. Given any set $X$, suppose we have a subset $\mathcal{S}\subseteq \mathcal{P}(X)$ of the powerset of $X$ that is closed under set difference $\setminus$. Then the structure $\langle \mathcal{S}; \setminus, \emptyset \rangle$ is a BCK-algebra. 

Another intended model, as perhaps indicated by the discussion in the introduction, is a collection of sentences together with $p\bdot q:= q\rightarrow p$ and $0:=\text{true}.$ In this context, the BCK-operation $\bdot$ is the order-dual of implication $\rightarrow$. Much of the early literature on BCK-algebras adopts the notational convention seen here using $\bdot$ and 0, though many later authors adopt the convention of using $\rightarrow$ and 1 instead. Iorgulescu calls our convention a ``right presentation'' and the presentation using $\rightarrow$ and 1 a ``left presentation'' in \cite{iorg07} and \cite{iorg08}. For this paper we will use Is\'{e}ki's original ``right presentation'' as given above.

We will need some other arithmetical properties of these algebras; the standard references for proofs of the following few results include the papers of Is\'{e}ki and Tanaka \cite{it78} and \cite{it78(2)} or the book \cite{mj94} by Meng and Jun. In an effort to keep this paper a bit more self-contained, we provide some proofs here. 

Given a BCK-algebra $\A$, define a relation $\leq$ on $\A$ by \[x\leq y \text{ if and only if } x\bdot y=0\,.\] We claim this is a partial order. Reflexivity follows from (BCK3) while antisymmetry follows from (BCK5). Transitivity requires a little more work.

\begin{proposition}[\hspace{1pt}\cite{it78}, Prop. 1] Let $\A$ be a BCK-algebra and $x,y,z\in\A$. Define $\leq$ as above.
\begin{enumerate}
\item If $x\leq y$, then $z\bdot y\leq z\bdot x$.
\item If $x\leq y$ and $y\leq z$, then $x\leq z$. Consequently, $\leq$ is a partial order.
\end{enumerate}
\end{proposition}

\begin{proof}\hfill
\begin{enumerate}\itemsep=.5em

\item Suppose $x\leq y$ so that $x\bdot y=0$. By (BCK1) we have
\[\bigl((z\bdot y)\bdot(z\bdot x)\bigr)\bdot 0=\bigl((z\bdot y)\bdot(z\bdot x)\bigr)\bdot(x\bdot y)=0\,.\]
But from (BCK4) we also have $0\bdot \bigl((z\bdot y)\bdot(z\bdot x)\bigr)=0$, and thus by (BCK5) it follows that 
\[(z\bdot y)\bdot(z\bdot x)=0\] which means $z\bdot y\leq z\bdot x$.

\item Suppose $x\leq y$ and $y\leq z$. Using (1) on the second inequality, we have $x\bdot z\leq x\bdot y=0$, meaning $(x\bdot z)\bdot 0=0$. But again, $0\bdot(x\bdot z)=0$ by (BCK4), and thus $x\bdot z=0$ by (BCK5). Hence, $x\leq z$.
\end{enumerate}

\end{proof}

Note that (BCK4) tells us that 0 is the least element of $\A$ with respect to the partial order $\leq$. The next result can be found in \cite[Thm 1]{it78} or the earlier paper \cite[p.19]{II66}, but the proof is too involved to include here.
\begin{theorem}\label{exchange}
In any BCK-algebra $\A$, we have \begin{align*}(x\bdot y)\bdot z=(x\bdot z)\bdot y\tag{ex}\end{align*} for all $x,y,z\in \A$.
\end{theorem}

\begin{lemma}[\hspace{1pt}\cite{it78}, Thm 2]\label{elem props} Let $\A$ be a BCK-algebra. Then
\begin{enumerate}
\item $x\bdot y\leq x$
\item $x\bdot 0=x$
\end{enumerate}
for all $x,y\in\A$.
\end{lemma}

\begin{proof}\hfill
\begin{enumerate}\itemsep=.5em
\item Using (ex), (BCK3), and (BCK4) we have
\[(x\bdot y)\bdot x=(x\bdot x)\bdot y =0\bdot y=0\] and hence $x\bdot y\leq x$.

\item Setting $y=0$ in (1) above yields $x\bdot 0\leq x$. On the other hand, from (BCK2) we have
\[\bigl(x\bdot (x\bdot 0)\bigr)\bdot 0=0\] and so $x\bdot (x\bdot 0)\leq 0$, from which it follows that $x\leq (x\bdot 0)$. Hence, $x\bdot 0=x$.
\end{enumerate}

\end{proof}

One may notice that, despite discussion of commutativity in the introduction, there is no mention of commutativity for the operation $\bdot$. In fact, the operation $\bdot$ is never commutative in a non-trivial BCK-algebra: if $x\neq 0$ then $x\bdot 0 =x\neq 0=0\bdot x$, for example.

In special cases, however, there is a term operation which \textit{is} commutative. Let $X$ be a set and consider the BCK-algebra $\langle \mathcal{P}(X); \setminus, \emptyset\rangle$. For $A,B\subseteq X$, it is easy to see that $A\setminus(A\setminus B) = A\cap B$ and thus $A\setminus(A\setminus B)=B\setminus(B\setminus A)$. Consider next the BCK-algebra of sentences with implication. Then $(p\rightarrow q)\rightarrow q=p\join q=(q\rightarrow p)\rightarrow p$. 

Of course, these two examples are particularly nice given their Boolean nature. But taking these examples as motivation, we define a term operation $\meet$ on a BCK-algebra $\A$ by $x\meet y:=x\bdot(x\bdot y)$. From Lemma \ref{elem props}(1) we see that $x\meet y\leq x$, and from (BCK2) we get $x\meet y\leq y$. Thus, the element $x\meet y$ is a common lower bound for $x$ and $y$, but in general is not the greatest lower bound of $x$ and $y$. For example, consider the algebra $\mb{M}_2$ whose Cayley table is given in Table \ref{tab:3}. The Hasse diagram for this algebra is the 3-element chain, $0<1<2$. So $\glb\{1,2\}=1$ while $2\meet 1 = 2\bdot (2\bdot 1)=0$. 
\begin{table}[h]
{\centering
\caption{\label{tab:3} The Cayley table for $\mb{M}_2$}
\begin{tabular}{c||ccc}
$\bdot$  & 0  & 1 & 2 \\\hline\hline
0           & 0  & 0 & 0 \\
1           & 1  & 0 & 0 \\
2           & 2  & 2 & 0
\end{tabular}\par}
\end{table} 

If $x\meet y=y\meet x$ we say that $x$ and $y$ {\bf commute}. If all pairs of elements of $\A$ commute then we say $\A$ is {\bf commutative}. 

\begin{remark}\label{comm elmts}
By (BCK3) and (2) of Lemma \ref{elem props} we have $x\meet x=x$ and $x\meet 0=0=0\meet x$ for all $x$ in any BCK-algebra $\A$.
\end{remark}

Recall that a class $\mathcal{K}$ of similar algebraic structures is a {\bf quasivariety} if $\mathcal{K}$ is closed under isomorphisms, subalgebras, direct products, and ultraproducts. A stronger notion is that of a variety: a class $\mathcal{K}$ of similar algebraic structures is a {\bf variety} if it is closed under homomorphic images, subalgebras, and direct products. It is well-known that any variety is a quasivariety. The class of all BCK-algebras is a quasivariety but not a variety \cite{wronski83}, while the subclass of commutative BCK-algebras is a variety \cite{yutani77}. One possible equational base which was given by Cornish in \cite[Thm 2.2]{cornish82} consists of (BCK1), (BCK4), $x\bdot 0=x$, and $x\bdot (x\bdot y)=y\bdot (y\bdot x)$. 

A function $f\colon \A\to \B$ between BCK-algebras is a {\bf BCK-homomorphism} if \[f(x\bdot_{\A} y)=f(x)\bdot_{\B}f(y)\,.\] Note that any BCK-homomorphism is 0-preserving by (BCK3). The {\bf kernel} of a BCK-homomorphism $f\colon \A\to\B$ is $\ker(f):=\{x\in\A\,\mid\, f(x)=0_{\B}\}$. We denote by $\mt{BCK}$ the category of BCK-algebras with BCK-homomorphisms while $\mt{cBCK}$ denotes the full subcategory of commutative BCK-algebras.

Let $\theta$ be a congruence on a BCK-algebra $\A$. The class $[0]_\theta$ of $\theta$ satisfies
\begin{enumerate}
\item[(i)] $0\in[0]_\theta$, and 
\item[(ii)] if $x\bdot y\in [0]_\theta$ and $y\in[0]_\theta$, then $x\in [0]_\theta$.
\end{enumerate}

\begin{proof}
The claim (i) is obvious. For (ii), suppose that $x\bdot y, y\in [0]_\theta$. Since $y\equiv_\theta 0$ and $\theta$ is a congruence, we have $x\bdot y\equiv_\theta x\bdot 0=x$. Then $x\bdot y\equiv_\theta 0$ implies that $x\equiv_\theta 0$.

\end{proof}

This motivates the following.
\begin{definition}
A subset $I\subseteq\A$ is an {\bf ideal} of $\A$ if $0\in I$ and the following implication is satisfied:
\begin{align*}\text{if $x\bdot y\in I$ and $y\in I$, then $x\in I$\,.}\tag{IP}\end{align*}
\end{definition} We note that the condition (IP), which stands for {\it ideal property}, is the order-dual of {\it modus ponens}. Let $\Id(\A)$ denote the set of ideals of $\A$.

\begin{lemma}
Given an ideal $I$ of $\A$, the relation $\theta_I$ on $\A$ defined by
\[x\,\theta_I\, y \text{ if and only if } x\bdot y, y\bdot x\in I\] is a congruence on $\A$ and $I=[0]_{\theta_I}$.
\end{lemma}

\begin{proof}
Symmetry of $\theta_I$ is clear. Reflexivity follows from that fact that $x\bdot x=0$ for all $x\in \A$ and $0\in I$. To see transitivity, suppose that $x\,\theta_I\, y$ and $y\,\theta_I\, z$. From (BCK1) we have
\[\bigl((x\bdot z)\bdot(x\bdot y)\bigr)\bdot(y\bdot z)=0\in I\,.\] Since $x\bdot y, y\bdot z\in I$, we apply (IP) twice to get $x\bdot z\in I$. An analogous argument shows $z\bdot x\in I$, and thus $x\,\theta_I\, z$. 

The compatibility of $\theta_I$ with the operation $\bdot$ uses (BCK1) together with other identities that will not be needed in this paper; we refer the reader to \cite[p.358]{it76}.

The last claim follows easily from the definition of $\theta_I$ and the identities $x\bdot 0=x$ and $0\bdot x=0$.

\end{proof}

On the other hand, in a BCK-algebra $\A$, a congruence $\theta$ is generally smaller than $\theta_{[0]_\theta}$; see \cite{wronski83} for an example. 

\begin{proposition}[\hspace{1pt}\cite{BR95}, Prop. 1] Let $\A$ be a BCK-algebra.
\begin{enumerate}
\item The set $\Id(\A)$ is a lattice with respect to subset inclusion, and the maps $I\mapsto \theta_I$ and $\theta\mapsto [0]_\theta$ are mutually inverse lattice isomorphisms between $\Id(\A)$ and the lattice $\Con_{\mt{BCK}}(\A)$ of relative BCK-congruences (congruences $\theta$ such that the quotient $\A/\theta$ is again a BCK-algebra).

\item $\A$ is relatively 0-regular.
\end{enumerate}
\end{proposition}

Given an ideal $I$ of $\A$, the quotient algebra $\A/\theta_I$ will be denoted by $\A/I$. Going forward, we adopt the somewhat unusual convention of writing $C_x$ to represent the element $[x]_I$ in $\A/I$. This notation is common in the older literature on BCK-algebras, and will make some computations involving pseudocommutators in quotient algebras less cluttered. Thus, the BCK-product in a quotient $\A/I$ is $C_x\bdot C_y=C_{x\bdot y}$. 

\begin{lemma}\label{ideals are downsets}
Let $I\in\Id(\A)$. Then $I$ is both a downset and a subalgebra.
\end{lemma}

\begin{proof} Take $x\in I$ and $a\in \A$ with $a\leq x$. Then $a\bdot x=0\in I$ and by (IP) we have $a\in I$. Thus, $I$ is a downset. Next, take $x,y\in I$. Since $x\bdot y\leq x$ and $I$ is a downset, we have that $x\bdot y\in I$, and so $I$ is also a subalgebra. 

\end{proof}

Example \ref{derived sub not ideal} gives an example of a subalgebra that is not an ideal.

We will write $\B\leq\A$ to indicate that $\B$ is a subalgebra of $\A$, and $I\unlhd\A$ to indicate that $I$ is an ideal of $\A$. Given a subset $S\subseteq \A$, the notations $\langle S\rangle$ and $(S]$ refer to the subalgebra generated by $S$ and the ideal generated by $S$, respectively. The following characterization of $(S]$ is Theorem 3 of \cite{it76}.
\begin{theorem}\label{ideal gen}
Given $S\subseteq \A$, we have
\begin{align*}
(S] = \{\,x\in \A\,\mid\, (\cdots((x\bdot s_1)\bdot s_2)\cdots )\bdot s_n=0 \text{ for some $s_1, \ldots, s_n\in S$}\,\}\,.
\end{align*}
\end{theorem}

\section{Pseudocommutators and the derived ideal}\label{commutators}

In this section, we use the notion of pseudocommutators to define the {\it derived ideal} of a BCK-algebra, which allows to define a functor which turns out to be left adjoint to the incluson functor $\mathsf{U}\colon \mt{cBCK}\to\mt{BCK}$.

Given categories $\mt{C}$ and $\mt{D}$, with $\mt{C}$ a full subcategory of $\mt{D}$, we say $\mt{C}$ is a {\bf reflective subcategory} of $\mt{D}$ if the inclusion functor $\mt{C}\to \mt{D}$ has a left adjoint; this adjoint functor is sometimes called a {\bf reflector}. 

Thus, this section culminates with a proof that $\mt{cBCK}$ is a reflective subcategory of $\mt{BCK}$.

\begin{remark} 
That $\mt{cBCK}$ is reflective in $\mt{BCK}$ can be proven by other means. These two categories are both {\it prevarieties} (closed under the formation of subalgebras and direct products) in the same signature, and $\mt{BCK}$ of course contains $\mt{cBCK}$. Then Theorem 3.6 of \cite{AS04} gives the result. In fact, equationally-definable subclasses of prevarieties are again prevarieties and thus reflective subcategories. What is new about the approach presented in this paper is the explicit construction of the reflector.

We thank Keith Kearnes for providing the argument above, as well as pointing the author to the reference \cite{AS04}.
\end{remark}

\begin{definition} The {\bf pseudocommutator} of $x$ and $y$ is 
\[[x,y]:=(x\meet y)\bdot(y\meet x)\,.\]
\end{definition}
This element measures how close $x$ and $y$ are to commuting, as we see in the next lemma.

\begin{lemma}\label{simple properties}
Let $\A$ be a BCK-algebra and $x,y\in\A$.
\begin{enumerate}
\item If $[x,y]=0$, then $x\meet y\leq y\meet x$.
\item If $x\meet y=y\meet x$, then $[x,y]=0$.
\item If $y\leq x$, then $[x,y]=0$.
\item $[x,y]\bdot x=0$ and $[x,y]\bdot y=0$.
\item $[x,0]=[0,x]=[x,x]=0$.
\end{enumerate}
\end{lemma} 

\begin{proof}\hfill
\begin{enumerate}

\item Apply Lemma \ref{elem props}(1).

\item Apply (BCK3).

\item Since $y\leq x$ we have $y\bdot x=0$, and so $y\meet x=y\bdot(y\bdot x)=y$. But $x\meet y\leq y$, so
\[[x,y]=(x\meet y)\bdot (y\meet x)=(x\meet y)\bdot y=0\,.\]

\item We use the identity (ex) together with the fact that $x\meet y\leq x$:
\[[x,y]\bdot x=\bigl((x\meet y)\bdot (y\meet x)\bigr)\bdot x=\bigl((x\meet y)\bdot x\bigr)\bdot (y\meet x)=0\bdot(y\meet x)=0\,.\] That $[x,y]\bdot y=0$ follows similarly.

\item This follows by Remark \ref{comm elmts} and (2) of this lemma.
\end{enumerate}

\end{proof}

Unfortunately, the inequality in Lemma \ref{simple properties}(1) cannot be upgraded to equality: looking ahead at Example \ref{derived sub not ideal}, we have $[3,2]=0$, but $3\meet 2=0<2=2\meet 3$. This is why we use the prefix ``pseudo.'' 

Given non-empty subsets $S,T\subseteq\A$, we define $[S,T]\subseteq\A$ by
\[[S,T]:=\Bigl\langle\,[s, t]\,\mid\, s\in S, t\in T\,\Bigr\rangle\,,\] the subalgebra generated by all pseudocommutators of elements of $S$ with elements of $T$. We will call the set of pseudocommutators $\{\,[s, t]\,\mid\, s\in S, t\in T\,\}$ a {\bf generating set} for $[S,T]$, and each pseudocommutator $[s,t]$ a {\bf generator} of $[S,T]$. 

The next lemma collects some expected properties of $[S,T]$ that will be used in later sections.

\begin{lemma}\label{derived products and inclusions}
\begin{enumerate}

\item Suppose $S\subseteq U\subseteq \A$ and $T\subseteq V\subseteq\A$ with $S$ and $T$ non-empty. Then $[S,T]\leq [U,V]$.

\item Let $S, U\subseteq \A$ and $T, V\subseteq \B$ with $S$, $U$, $T$, and $V$ non-empty. Then 
\[[S\times T, U\times V]_{\A\times\B} \leq [S,U]_\A\times [T,V]_\B\,.\]
\end{enumerate}
\end{lemma}

\begin{proof}\hfill

(1) Any generator $[s,t]$ of $[S,T]$ is also a generator of $[U,V]$. The claim follows.

(2) Take $(s, t)\in S\times T$ and $(u,v)\in U\times V$. Then
\begin{align*}
\bigl[\,(s,t)\,,\,(u,v)\,\bigr]_{\A\times\B} &= \bigl(\,(s,t)\meet (u,v)\,\bigr)\bdot\bigl(\,(u,v)\meet (s,t)\,\bigr)\\
&=\bigl(s\meet u\,,\, t\meet v\bigr)\bdot \bigl(u\meet s\,,\, v\meet t\bigr)\\
&=\Bigl(\,(s\meet u)\bdot_\A(u\meet s)\,,\,(t\meet v)\bdot_\B(v\meet t)\,\Bigr)\\
&=\Bigl(\,[s,u]_\A\,,\, [t,v]_\B\,\Bigr)\,.
\end{align*} This implies that $[S\times T, U\times V]_{\A\times\B}\leq [S,U]_\A\times [T,V]_\B$, as claimed.

\end{proof}

The {\bf derived subalgebra} of $\A$ is $\A':=\A^{(1)}:=[\A, \A]$. From there we define the {\bf derived series} of $\A$:
\[\A^{(2)}=[\A', \A']\text{ and } \A^{(k+1)}=[\A^{(k)}, \A^{(k)}]\]
for $k>2$. It is easy to see that $\A^{(k+1)}$ is a subalgebra of $\A^{(k)}$, so this is a descending series of subalgebras:
\[\A\geq \A'\geq \A^{(2)}\geq \cdots\geq \A^{(k)}\geq\cdots\,.\]
We say $\A$ is {\bf solvable} if $\A^{(n)}=\{0\}$ for some $n\in\bb{N}$.

Najafi proved in \cite{najafi13} that $\A$ is commutative if and only if $\A'=\{0\}$, so any commutative BCK-algebra is solvable. In the subsequent paper \cite{NS14}, Najafi and Saied proved that a quotient $\A/I$ is commutative if and only if $\A'\leq I$. The proof of this follows from the fact that $[C_x, C_y]_{\A/I}= C_{[x,y]_\A}$.

In general, the derived subalgebra $\A'$ need not be an ideal, see Example \ref{derived sub not ideal}. To remedy this, we introduce the following:
\begin{definition}
The {\bf derived ideal} of $\A$ is the ideal generated by $\A'$. Symbolically,
\[\DI(\A):= \bigl(\A'\bigr] = \bigcap_{\substack{\A'\subseteq I\\I\unlhd\A}}I\,.\]
\end{definition} From Theorem \ref{ideal gen}, it follows that \[\DI(\A)=\{\,x\in \A\,\mid\, (\cdots((x\bdot a_1)\bdot a_2)\cdots )\bdot a_n=0 \text{ for some $a_1, \ldots, a_n\in \A'$}\,\}\,.\]
Notice that $\A'\leq \DI(\A)$ and thus the quotient $\A/\DI(\A)$ is commutative.

\begin{example}\label{derived sub not ideal}
Consider the algebra $\A$ whose Cayley table and $\meet$-table are shown in Table \ref{tab:alg1}. The Hasse diagram is shown in Figure \ref{fig:A}. From the $\meet$-table, one computes $\A'=\{0, 1, 2\}$ which is a subalgebra but not an ideal: notice that $4\bdot 2=2\in \A'$ and $2\in\A'$, but $4\notin\A'$, meaning $\A'$ does not satisfy (IP). The derived ideal is then $\DI(\A)=\{0, 1, 2, 4\}$ and the quotient $\A/\DI(\A)$ is isomorphic to the (unique) two-element BCK-algebra, which is commutative. 
\begin{table}[h]
{\centering
\caption{\label{tab:alg1} The Cayley table and $\meet$-table for $\A$}
\begin{tabular}{c||ccccc}
$\bdot$  & 0  & 1 & 2 & 3 & 4\\\hline\hline
0           & 0  & 0 & 0 & 0 & 0\\
1           & 1  & 0 & 1 & 0 & 1\\
2           & 2  & 2 & 0 & 0 & 0\\
3		& 3 & 3 & 3 & 0 & 3 \\
4		& 4 & 4 & 2 & 2 & 0
\end{tabular}
\hspace{1cm}
\begin{tabular}{c||ccccc}
$\meet$  & 0  & 1 & 2 & 3 & 4\\\hline\hline
0           & 0  & 0 & 0 & 0 & 0\\
1           & 0  & 1 & 0 & 1 & 0\\
2           & 0  & 0 & 2 & 2 & 2\\
3		& 0 & 0 & 0 & 3 & 0\\
4		& 0 & 0 & 2 & 2 & 4
\end{tabular}\par
}
\end{table}
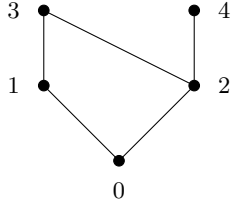
\begin{figure}[h]
\centering
\begin{tikzpicture}
\filldraw (0,0) circle (2pt);
\filldraw (-1,1) circle (2pt);
\filldraw (1,1) circle (2pt);
\filldraw (-1,2) circle (2pt);
\filldraw (1,2) circle (2pt);
\draw [-] (0,0) -- (-1,1) -- (-1,2);
\draw [-] (0,0) -- (1,1) -- (1,2);
\draw [-] (1,1) -- (-1,2);
	\node at (0,-.4) {\small 0};
	\node at (-1.4, 1) {\small 1};
	\node at (1.4, 1) {\small 2};
	\node at (-1.4, 2) {\small 3};
	\node at (1.4, 2) {\small 4};
\end{tikzpicture}
\caption{Hasse diagram for $\A$}
\label{fig:A}
\end{figure}

\end{example}
 
\begin{definition}
The {\bf commutativization} of $\A$ is the quotient algebra \[\A^{\comm}:=\A/\DI(\A)\,.\]
\end{definition}

The commutativization enjoys the following universal property:

\begin{theorem}\label{universal prop}
Let $\A\in\mt{BCK}$. Given any homomorphism $\phi\colon \A\to \C$ with $\C\in\mt{cBCK}$, there exists a unique homomorphism $\psi\colon \A^{\comm}\to \C$ such that $\phi=\psi\circ \pi$, where $\pi$ is the natural projection $\A\to\A^{\comm}$.
\end{theorem}

\begin{center}
\begin{tikzcd}
\A \arrow[rr, "\phi"] \arrow[dd, "\pi"] & &\C\\
&&\\
\A^{\comm} \arrow[rruu, "{\exists!\, \psi}", dashed]
\end{tikzcd}
\end{center}

\begin{proof}
Since $\phi$ is a homomorphism it preserves BCK-terms. Together with the fact that $\C$ is commutative, we see
\[\phi\bigl([x,y]_\A\bigr)=\bigl[\phi(x), \phi(y)\bigr]_\C=0_\C\]
and thus $[x,y]_\A\in\ker(\phi)$ for all $x,y\in \A$. This means $\A'\subseteq \ker(\phi)$, but since $\ker(\phi)$ is an ideal we get $\ker(\pi)=\DI(\A)\subseteq \ker(\phi)$ as well. By Theorem 15 of \cite{it78}, there exists a unique homomorphism $\psi\colon \A^{\comm}\to \C$ such that $\phi=\psi\circ \pi$.

\end{proof}

\begin{theorem}\label{hom props of derived ideal}
Let $\phi\colon\A\to\B$ be a homomorphism of BCK-algebras. Then:
\begin{enumerate}
\item $\phi\bigl(\DI(\A)\bigr)\subseteq \DI(\B)$, and 
\item $\phi$ induces a homomorphism $\Phi\colon \A^{\comm}\to\B^{\comm}$ by 
\[\Phi(C_x)=C_{\phi(x)}\,.\]
\item Consequently, $\mt{cBCK}$ is a reflective subcategory of $\mt{BCK}$ with the commutativization functor as the reflector.
\end{enumerate}
\end{theorem}

\begin{proof}\hfill

(1) We first show that $\phi(\A')\subseteq \B'$. Take $y\in \phi(\A')$. Then there exists $x\in \A'$ with $\phi(x)=y$. Since $x\in \A'$, there exists an $n\in\bb{N}$ and a BCK-term $u$ with $n$ variables such that we can write 
\[x=u^\A\bigl(\,[a_1,c_1],\ldots, [a_n,c_n]\,\bigr)\] for some $a_1,\ldots, a_n,c_1,\ldots, c_n\in\A$. 
But then 
\begin{align*}
y=\phi(x)&=\phi\Bigl(u^\A\bigl(\,[a_1,c_1],\ldots, [a_n,c_n]\,\bigr)\Bigr)\\
&=u^\B\Bigl(\phi\bigl([a_1,c_1]\bigr), \ldots, \phi\bigl([a_n,c_n]\bigr)\Bigr)\\
&=u^\B\Bigl(\bigl[\phi(a_1),\phi(c_1)\bigr],\ldots,\bigl[\phi(a_n),\phi(c_n)\bigr]\Bigr)\in\B'\,.
\end{align*}
Hence, $\phi(\A')\subseteq \B'$.

Now take $y\in \phi\bigl(\DI(\A)\bigr)$. Then there is $x\in \DI(\A)$ with $\phi(x)=y$, and since $x\in \DI(\A)$ there exist $a_1, \ldots, a_n\in\A'$ such that
\[(\cdots ((x\bdot a_1)\bdot a_2)\cdots)\bdot a_n=0_\A\,.\] Applying $\phi$ to both sides gives
\[y\bdot \phi(a_1)\bdot \phi(a_2)\bdot\cdots \bdot \phi(a_n) =0_\B\,,\] where the product is left-associated. Since each $\phi(a_i)\in \B'$, we have $y\in DI(\B)$. Thus, $\phi\bigl(\DI(\A)\bigr)\subseteq \DI(\B)$.

(2) We show that $\Phi$ is well-defined. Suppose $C_x=C_y$ for some $x,y\in \A$. This means that $x\bdot y, y\bdot x\in \DI(\A)$. But then by (1) above, we have $\phi(x)\bdot\phi(y), \phi(y)\bdot\phi(x)\in\DI(\B)$. Hence, $C_{\phi(x)}=C_{\phi(y)}$. An easy computation now shows that $\Phi$ is a BCK-homomorphism.

(3) It is straightforward to see that $(-)^{\comm}$ is a covariant functor $\mt{BCK}\to\mt{cBCK}$ which assigns to each BCK-algebra $\A$ its commutativization $\A^{\comm}$ and to each homomorphism $\phi\colon \A\to\B$ the induced map $\Phi\colon \A^{\comm}\to\B^{\comm}$. From Theorem \ref{universal prop} we see that this functor is left adjoint to the inclusion functor $\mathsf{U}\colon \mt{cBCK}\to\mt{BCK}$. Hence, $\mt{cBCK}$ is a reflective subcategory of $\mt{BCK}$ with reflector $(-)^{\comm}$.

\end{proof}

\section{Central series and nilpotence}\label{nilpotence}

In this section, we introduce the concept of central series for BCK-algebras and use them to define a notion of nilpotence, which is then connected to pseudocommutators. As mentioned in the introduction, we view nilpotence as a measure of how commutative a BCK-algebra is: between two algebras $\A$ and $\B$, the algebra with lower nilpotence class is considered ``more commutative.''

\begin{definition} Let $\A$ be a BCK-algebra. 
\begin{enumerate}
\item An {\bf ascending central series} of $\A$ is a sequence of subalgebras
\[\{0\}=\B_0\leq \B_1\leq \cdots \leq \B_{n-1}\leq \B_n\leq\cdots\leq \A \] such that $[\B_{i+1}, \A]\leq \B_i$. 

\item A {\bf descending central series} of $\A$ is a sequence of subalgebras
\[\{0\}\leq \cdots \leq \B_n\leq \B_{n-1}\leq \cdots \leq \B_1\leq \B_0=\A \] such that $[\B_{i}, \A]\leq \B_{i+1}$.

\item We say $\A$ is {\bf nilpotent} if it has a central series of finite length. The length of a finite central series is the {\bf nilpotence class} of $\A$. 
\end{enumerate}
\end{definition}

\begin{remark} Readers familiar with universal algebra may know that there is a notion of {\it nilpotence} for general algebraic structures which is defined using commutators of congruences; the precise definition and accompanying details can be found in Ch.\ 6 of Freese and McKenzie's text \cite{FM87}. However, BCK-algebras are congruence distributive algebras, see Theorem 1 of \cite{palasinski81}; so the commutator of two congruences is just their intersection. This implies that no non-trivial BCK-algebra is nilpotent under Freese and McKenzie's definition.
\end{remark}

Just as in group theory, there are other equivalent definitions of nilpotence. First, for $c\geq 3$, the {\bf (left-normed) pseudocommutator of weight $c$} is defined by
\[[x_1, x_2, \ldots, x_c]=\bigl[\,[x_1,\ldots, x_{c-1}]\,,\, x_c\,\bigr]\,.\] Using these higher order pseudocommutators, we define the {\bf upper central series} of $\A$ by $Z_0(\A)=\{0\}$, 
\[Z_1(\A) = \langle\, x\in\A \,\mid\,[x,y]=0\text{ for all $y\in\A$}\,\rangle\,,\]
and for $n>1$,
\[Z_n(\A) = \langle\, x\in\A \,\mid\,[x,y_1,\ldots, y_n]=0\text{ for all $y_1,\ldots, y_n\in\A$}\,\rangle\,.\] When the underlying algebra is understood, we will simply write $Z_k$ for $Z_k(\A)$.

\begin{lemma}
The upper central series is an ascending central series.
\end{lemma}

\begin{proof}
It is straightforward to show that this series of subalgebras is ascending. We show $[Z_k, \A]\leq Z_{k-1}$. 

Take a generator $z\in [Z_k, \A]$ so that $z=[x,y]$ where $x\in Z_k$ and $y\in\A$. Now take $y_1, \ldots, y_{k-1}\in \A$. Then we have
\[[z, y_1, \ldots, y_{k-1}]=\bigl[[x,y],y_1,\ldots, y_{k-1}\bigr]=[x, y, y_1,\ldots, y_{k-1}]=0\] since $x\in Z_k$, and this equality tells us $z\in Z_{k-1}$. Since all the generators of $[Z_k,\A]$ are elements of $Z_{k-1}$, we have $[Z_k, \A]\leq Z_{k-1}$.

\end{proof}

We call $Z_1(\A)$ the {\bf pseudocenter} of $\A$, as the elements of $Z_1$ need not be central, as the next example shows.

\begin{example}\label{big example} Consider the algebra $\B$ defined by Table \ref{tab:alg3}. Some computation reveals that $Z_1(\B)=\{0, 1, 3,4,5,7\}$ while the central elements of $\B$ are $\{0, 1, 4,5\}$.

This example also shows that $Z_1$ need not be an ideal: note that $2\bdot 1=1\in Z_1(\B)$ and $1\in Z_1(\B)$, but $2\notin Z_1(\B)$.
\begin{table}[h]
{\centering
\caption{\label{tab:alg3} Cayley table and $\meet$-table for $\B$}
\begin{tabular}{c||cccccccc}
$\bdot$  & 0  & 1 & 2 & 3 & 4 & 5 & 6 & 7\\\hline\hline
0           & 0  & 0 & 0 & 0 & 0 & 0 & 0 & 0\\
1           & 1  & 0 & 0 & 0 & 1 & 0 & 0 & 0\\
2           & 2  & 1 & 0 & 0 & 2 & 1 & 0 & 0\\
3		& 3  & 1 & 1 & 0 & 3 & 1 & 1 & 0\\
4		& 4  & 4 & 4 & 4 & 0 & 0 & 0 & 0\\
5		& 5  & 4 & 4 & 4 & 1 & 0 & 0 & 0\\
6		& 6  & 5 & 4 & 4 & 2 & 1 & 0 & 0\\
7		& 7  & 5 & 5 & 4 & 3 & 1 & 1 & 0\\
\end{tabular}
\hspace{1cm}
\begin{tabular}{c||cccccccc}
$\meet$  & 0  & 1 & 2 & 3 & 4 & 5 & 6 & 7\\\hline\hline
0           & 0  & 0 & 0 & 0 & 0 & 0 & 0 & 0\\
1           & 0  & 1 & 1 & 1 & 0 & 1 & 1 & 1\\
2           & 0  & 1 & 2 & 2 & 0 & 1 & 2 & 2\\
3		& 0  & 1 & 1 & 3 & 0 & 1 & 1 & 3\\
4		& 0  & 0 & 0 & 0 & 4 & 4 & 4 & 4\\
5		& 0  & 1 & 1 & 1 & 4 & 5 & 5 & 5\\
6		& 0  & 1 & 2 & 2 & 4 & 5 & 6 & 6\\
7		& 0  & 1 & 1 & 3 & 4 & 5 & 5 & 7\\
\end{tabular}\par
}
\end{table}
\end{example}

A natural example of a descending central series is the {\bf lower central series} defined by $\A_{(1)}=\A'$, $\A_{(2)}=[\A', \A]$, and in general, $\A_{(k+1)}=[\A_{(k)}, \A]$ for $k>2$. It is straightforward to check that $\A_{(k+1)}$ is a subalgebra of $\A_{(k)}$, and it is a central series by definition.

We will show that an algebra $\A$ is nilpotent if and only if the lower central series is finite if and only if the upper central series is finite, and that the two series will have the same length.

\begin{theorem}\label{central series props}
\begin{enumerate}
\item  Suppose $\{\B_n\}$ is a descending central series of $\A$. Then $\A_{(k)}\leq \B_k$ for all $k$. That is, the lower central series is the fastest descending central series.

\item Suppose $\{\B_n\}$ is an ascending central series of $\A$. Then $\B_k\leq Z_k(\A)$ for all $k$. That is, the upper central series is the fastest ascending central series. 

\item $Z_n(\A)=\A$ for some $n\in\bb{N}$ if and only if $\A_{(n)}=\{0\}$, in which case the least $n$'s for which these occur are the same.
\end{enumerate}
\end{theorem}

\begin{proof}\hfill

(1) Note that $\A_{(1)}=[\A,\A]=[\B_0, \A]\leq \B_1$, giving our base case. Suppose that $\A_{(k-1)}\leq \B_{k-1}$ for some $k\geq 2$. Then by Lemma \ref{derived products and inclusions} we have
\[\A_{(k)}=[\A_{(k-1)}, \A]\leq [\B_{k-1}, \A]\leq \B_k\] and the result follows by induction.

(2) By definition we have $\B_0=\{0\}=Z_0$. Suppose $\B_{k-1}\leq Z_{k-1}$ for some $k\geq 1$. Take $x\in \B_k$ and $y_1, y_2, \ldots, y_k\in\A$. Then $[x,y_1]\in [\B_k,\A]\leq \B_{k-1}\leq Z_{k-1}$ and we have
\[[x,y_1,\ldots, y_k]=\bigl[[x, y_1],y_2,\ldots, y_k\bigr]=0\,,\] which implies $x\in Z_k$. Thus $\B_k\leq Z_k$ and the result follows by induction.

(3) Suppose $Z_n=\A$. By reversing indices, we can view the upper central series as a descending central series: set $\B_k=Z_{n-k}$ for $k=0,1, \ldots, n$. By (1), we get $\A_{(k)}\leq \B_k=Z_{n-k}$, which means $\A_{(n)}=\{0\}$. 

The converse is analogous using (2).

\end{proof}

\begin{corollary}
A BCK-algebra $\A$ is nilpotent if and only if the upper central series terminates at $\A$ in finitely steps if and only if the lower central series terminates at $\{0\}$ in finitely many steps.
\end{corollary}

\begin{corollary}
A BCK-algebra $\A$ is commutative if and only if $Z_1(\A)=\A$. That is, the commutative BCK-algebras are precisely those that are nilpotent of class 1.
\end{corollary}

So we see that any commutative BCK-algebra is nilpotent, but the next example shows a nilpotent BCK-algebra need not be commutative.

\begin{example}\label{big example2}
Recall the algebra $\B$ from Example \ref{big example}. We saw previously that $Z_1(\B)=\{0, 1, 3,4,5,7\}$. Some computation reveals that
\[\B_{(1)}=\{0,1\}\subset Z_1(\B)\] and further that $Z_2(\B)=\B$ and $\B_{(2)}=\{0\}$. Thus, $\B$ is nilpotent of class 2, using any definition. However, $\B$ is not commutative since $3\meet 6 = 1\neq 2=6\meet 3$.

\end{example}

\begin{example}\label{Mn's}
Consider the algebra $\M_3$ whose Cayley table is shown in Table \ref{tab:M3}.  
\begin{table}[h]
{\centering
\caption{\label{tab:M3} The Cayley table and $\meet$-table for $\M_3$}
\begin{tabular}{c||cccc}
$\bdot$  & 0  & 1 & 2 & 3\\\hline\hline
0           & 0  & 0 & 0 & 0\\
1           & 1  & 0 & 0 & 0\\
2           & 2  & 2 & 0 & 0\\
3		& 3 & 3 & 3 & 0
\end{tabular}
\hspace{1cm}
\begin{tabular}{c||cccc}
$\meet$  & 0  & 1 & 2 & 3\\\hline\hline
0           & 0  & 0 & 0 & 0\\
1           & 0  & 1 & 1 & 1\\
2           & 0  & 0 & 2 & 2\\
3		& 0 & 0 & 0 & 3
\end{tabular}\par
}
\end{table} 

From this we compute
\begin{align*}
\M_3'&=\{0,1,2\}\\
(\M_3)_{(2)}&=\{0,1\}\\
(\M_3)_{(3)}&=\{0\}
\end{align*} so $\M_3$ is nilpotent of class 3. We note that the Hasse diagram for $\M_3$ is a four-element chain.

This example can be generalized as follows: let $(M_n,\leq)$ be the $(n+1)$-element chain with carrier set $\{0,1,\ldots, n\}$, where $\leq$ is the natural order. If we define an operation $\bdot$ on $M_n$ by
\begin{align*}
x\bdot y=\begin{cases} 0 & \text{ if $x\leq y$}\\x & \text{ otherwise}\end{cases}\,,\tag{$\dagger$}
\end{align*} then $\M_n:=\langle M_n; \bdot, 0\rangle$ is a BCK-algebra such that the order induced by $\bdot$ coincides with $\leq$. (That these are BCK-algebras follows from Theorem 4 of \cite{iseki75}: the algebra $\M_n$ is obtained from $\M_{n-1}$ by appending a new top element and suitably extending the operation. We note here also that any poset $(P,\leq)$ with a least element, when equipped with the operation defined by ($\dagger$), is a BCK-algebra whose induced order coincides with $\leq$; this is mentioned in \cite[p. 29]{GL92}, though the authors leave the proof as an exercise to the reader.)

One can compute that $(\M_n)_{(k)}=\{0,1,\ldots, n-k\}$ for $1\leq k\leq n$, and thus $(\M_n)_{(n)}=\{0\}$. So $\M_n$ is nilpotent of class $n$. This shows there are BCK-algebras with arbitrarily large nilpotence class.
\end{example}

\begin{remark} Recall the notion of {\it commuting degree} from the introduction. For a finite BCK-algebra, this is the quotient \[\frac{\text{(\# pairs of commuting elements)}}{\text{(total \# of pairs of elements)}}\,,\] which is the probability that two randomly selected elements commute. In \cite{evans25} it is proven that the algebra $\M_n$ is the only BCK-algebra of order $n+1$ (up to isomorphism) that achieves the minimum possible commuting degree. That is, among algebras of order $n+1$, these algebras are ``as far as possible'' from being commutative. Therefore we think it is interesting to note that the algebra $\M_n$ also has the largest possible nilpotence class among algebras of order $n+1$, reinforcing the observation that this algebra is far from commutative. 

That being said, while this suggests a connection between commuting degree and nilpotence class, a close study of the possible relationship between these two quantities has not yet been undertaken.

\end{remark}

\begin{example}\label{unit interval}
On the other hand, consider the unit interval $I=[0,1]\subseteq \bb{R}$. By our comments in Example \ref{Mn's}, $I$ becomes a BCK-algebra when equipped with the operation defined by ($\dagger$); we denote this algebra by $\I$. According to \cite[Ex. 8]{NS14}, the algebra $\I$ is not solvable, and thus not nilpotent by Proposition \ref{nil implies solvable}. But this is also easy to see directly from the definition: one computes 
\[[x,y]=\begin{cases}x &\text{ if $x\leq y$}\\0 & \text{ otherwise}\end{cases}\] from which it follows that $\I_{(n)}=[0,1)$ for all $n\in\bb{N}$. Thus, $\I$ is not nilpotent. We note that $\I$ contains, for each $n$, a subalgebra isomorphic to $\M_n$.
\end{example}

\begin{proposition}\label{nil implies solvable}
If $\A$ is a nilpotent BCK-algebra, then it is also solvable.
\end{proposition}

\begin{proof}
We claim that $\A^{(n)}\leq \A_{(n)}$ for all $n\in\bb{N}$. The base case is obvious since $\A^{(1)}:=\A':=\A_{(1)}$.

Suppose $\A^{(k)}\leq \A_{(k)}$ for some $k>1$. By the definitions and an application of Lemma \ref{derived products and inclusions}, we have
\[\A^{(k+1)}=[\A^{(k)}, \A^{(k)}]\leq [\A_{(k)}, \A]=\A_{(k+1)}\,.\] The claim follows by induction, and from this it is immediate that if $\A$ is nilpotent, then $\A$ is solvable.

\end{proof}

Thus, we have the following containments
\[\text{Commutative BCK }\subset\text{ Nilpotent BCK }\subseteq
\text{ Solvable BCK }\subset\text{ BCK }\,.\]

We have not yet found an example of a solvable BCK-algebra that is not nilpotent, though we suspect such examples exist. Any such example must be infinite (see Theorem \ref{fin = nil} below) and non-commutative.

\section{Properties of nilpotent algebras}\label{nilp props}

In this section we discuss some of the varietal properties of nilpotent BCK-algebras. We close the paper by showing that every BCK-algebra of finite height is nilpotent.

\begin{theorem}\label{not quite pseudo}
The class of nilpotent BCK-algebras is closed under subalgebras, finite products, and surjective homomorphisms.
\end{theorem}

\begin{proof}
Let $\A$ be a BCK-algebra and suppose $\B\leq \A$. We prove that $\B_{(n)}\leq \A_{(n)}$ for all $n\in\bb{N}$, which gives the first claim of the theorem. 

It is clear that $\B'\leq \A'$ since any generator of $\B'$ is a generator of $\A'$. So now suppose $\B_{(k-1)}\leq\A_{(k-1)}$ for some $k>1$. Take a generator $z\in \B_{(k)}$. Then we can write $z=[x,y]$ for some $x\in \B_{(k-1)}$ and $y\in\B$. By the inductive hypothesis, $x\in \A_{k-1}$ and of course $y\in \A$. Thus, $z\in \A_{(k)}=[\A_{(k-1)}, \A]$. Since all generators of $\B_{(k)}$ are contained in $\A_{(k)}$, it follows that $\B_{(k)}\leq\A_{(k)}$. By induction we have $\B_{(n)}\leq \A_{(n)}$, and from this it follows that if $\A$ is nilpotent, then so is $\B$.

Next, we consider finite products. Let $\A$ and $\B$ be BCK-algebras. Applying Lemma \ref{derived products and inclusions}(2), we see
\[(\A\times \B)'=[\A\times \B, \A\times \B]_{\A\times \B}\leq [\A,\A]_{\A}\times [\B, \B]_{\B}=\A'\times\B'\,.\]
Then, using Lemma \ref{derived products and inclusions} together with the observation above, we have
\begin{align*}
(\A\times \B)_{(2)}
=[(\A\times\B)', \A\times\B]_{\A\times\B}
&\leq [\A'\times\B', \A\times\B]_{\A\times\B}\\
&\leq [\A', \A]_{\A}\times[\B', \B]_{\B}\\
&=\A_{(2)}\times\B_{(2)}\,.\end{align*}
The same argument can be applied inductively to see that $(\A\times \B)_{(n)}\leq \A_{(n)}\times \B_{(n)}$. From this it follows that if $\A$ and $\B$ are both nilpotent, then so is $\A\times \B$. We extend this to more factors by another induction. 

From the above paragraph we see that if $\A_1,\ldots, \A_k$ are nilpotent of classes $c_1, \ldots, c_k$, respectively, then $\prod_{i=1}^k\A_i$ is nilpotent of class $\max\{c_1,\ldots,c_k\}$.

Finally, we show that if $\A$ is nilpotent and $\phi\colon \A\to\B$ is a surjective homomorphism, then $\B$ is nilpotent. To this end, we first prove that $\phi(\A_{(n)})=\B_{(n)}$ for all $n\in\bb{N}$.

We saw in the proof of Lemma \ref{hom props of derived ideal}(1) that $\phi(\A')\subseteq \B'$. For the other inclusion, take $b\in \B'$. Then we can write 
\[b=u^\B\bigl([y_1,z_1],\ldots,[y_n,z_n]\bigr)\] for some BCK-term $u$ with $n$ variables and elements $y_1,\ldots, y_n,z_1,\ldots, z_n\in \B$. Since $\phi$ is surjective, there exist $a_i, c_i\in \A$ such that $y_i=\phi(a_i)$ and $z_i=\phi(c_i)$ for $i=1,\ldots, n$, and we have
\begin{align*}
b=u^\B\bigl([y_1,z_1],\ldots,[y_n,z_n]\bigr)
&=u^\B\bigl([\phi(a_1),\phi(c_1)],\ldots,[\phi(a_n),\phi(c_n)]\bigr)\\
&=u^\B\Bigl(\phi\bigl([a_1,c_1]\bigr),\ldots,\phi\bigl([a_n,c_n]\bigr)\Bigr)\\
&=\phi\Bigl(\,u^\A\bigl([a_1,c_1],\ldots,[a_n,c_n]\bigr)\,\Bigr)\in\phi(\A')\,.
\end{align*} Thus, $\B'\subseteq \phi(\A')$ and this gives $\B'=\phi(\A')$. The rest of the argument is similar and proceeds by an induction. Therefore $\phi(\A_{(n)})=\B_{(n)}$ for all $n\in\bb{N}$. 

Now if $\A$ is nilpotent then $\A_{(k)}$ for some $k$, and so $\phi(\A_{(k)})=\phi(\{0_\A\})=\{0_\B\}=\B_{(k)}$. So $\B$ is nilpotent as well.

\end{proof}

We recall that a class of algebras $\mc{K}$ is called a
\begin{itemize}
\item {\bf variety} if $\mc{K}$ is closed under homomorphic images, subalgebras, and direct products.

\item {\bf quasivariety} if $\mc{K}$ is closed under isomorphisms, subalgebras, direct products, and ultraproducts.

\item {\bf pseudovariety} if $\mc{K}$ is closed under homomorphic images, subalgebras, and finitary direct products.
\end{itemize}

Example \ref{infinite product} below shows that an infinite product of nilpotent BCK-algebras need not be nilpotent, so the class of nilpotent BCK-algebras, call it $\mt{BCK}_{\nil}$, is neither a variety nor a quasivariety. Theorem \ref{not quite pseudo} shows that $\mt{BCK}_{\nil}$ is {\it almost} -- but not quite -- a pseudovariety. As mentioned in Section \ref{prelims}, the class $\mt{BCK}$ is not a variety: specifically, Wro\'{n}ski showed in \cite{wronski83} that $\mt{BCK}$ is not closed under homomorphic images. However, Blok and Raftery showed in \cite{BR95} that the lattice of ideals of a BCK-algebra is isomorphic to the lattice of BCK-congruences (that is, congruences $\theta$ of $\A$ such that the quotient $\A/\theta$ is again a BCK-algebra). Hence, by the previous theorem, if $\A$ is nilpotent then so is $\A/I$ for any ideal of $\A$. 

Thus, if $\mt{V}$ is a variety of BCK-algebras, then the subclass $\mt{V}_{\nil}$ of nilpotent members of $\mt{V}$ is itself closed under homomorphic images, and hence $\mt{V}_{\nil}$ is a pseudovariety. Still, though, Example \ref{infinite product} shows that $\mt{V}_{\nil}$ need not be a variety (or quasivariety) in general. In particular, any variety $\mt{V}$ of BCK-algebras that contains the algebra $\I$ of Example \ref{unit interval} will contain the algebra $\D$ defined below, and thus $\mt{V}_{\nil}$ will not be a variety (or quasivariety).

\begin{example}\label{infinite product}
Recall the algebras $\M_n$ defined in Example \ref{Mn's} and consider $\D=\prod_{i=3}^\infty\M_i$. Label a sequence of elements in $\D$:
\begin{align*}
\alpha_0&=(3,4,5,6\ldots)\\
\alpha_1&=(2,3,4,5\ldots)\\
\alpha_2&=(1,2,3,4\ldots)\\
\alpha_3&=(0,1,2,3\ldots)\\
\alpha_4&=(0,0,1,2,\ldots)\\
&\vdots
\end{align*}
One computes
\begin{align*}
[\alpha_1, \alpha_0]&=(\alpha_1\meet\alpha_0)\bdot(\alpha_0\meet\alpha_1)=\alpha_1\bdot\mathbf{0}=\alpha_1\in \D_{(1)}\\
[\alpha_2, \alpha_1]&=(\alpha_2\meet\alpha_1)\bdot(\alpha_1\meet\alpha_2)=\alpha_2\bdot\mathbf{0}=\alpha_2\in \D_{(2)}
\end{align*} and in general
\begin{align*}
[\alpha_{k+1}, \alpha_k]=\alpha_k\in \D_{(k)}\,.
\end{align*} Thus, $\D_{(n)}\neq\{\mathbf{0}\}$ for all $n\in\bb{N}$, meaning $\D$ is not nilpotent.
\end{example}

One may recall from group theory that the class of all nilpotent groups is not a variety, but the class of all groups with nilpotence class at most $c$ {\it is} a variety. The situation for BCK-algebras is similar: as discussed previously, $\mt{BCK}_{\nil}$ is not a quasivariety, but the class $\mt{BCK}_c$ of nilpotent BCK-algebras with nilpotence class at most $c$ {\it is} a quasivariety.

\begin{lemma}
$\A\in\mt{BCK}_{c}$ if and only if $[x_1, x_2, \ldots, x_c]=0$ for all choices of $x_1, \ldots, x_c\in \A$.
\end{lemma}

\begin{proof}
Fix $c\in\bb{N}$. Suppose $[a_1, a_2, \ldots, a_c]=0$ for all $a_1, \ldots, a_c\in \A$. This immediately implies that $Z_c(\A)=\A$ and thus $\A\in\mt{BCK}_{c}$

Conversely, suppose $[a_1, a_2, \ldots, a_c]\neq 0$ for some $a_1,\ldots, a_c\in \A$. This tells us $\A_c\neq \{0\}$. If $[a_1, \ldots, a_i]=0$ for any $1\leq i\leq c$, then we would have
\[[a_1, \ldots, a_c]=\bigl[[a_1,\ldots, a_i], a_{i+1},\ldots , a_c\bigr]=[0,a_{i+1},\ldots, a_c]=0\,,\] a contradiction. So $[a_1, \ldots, a_i]\neq 0$ for each $1\leq i\leq c$, which means $\A_i\neq\{0\}$ for each $1\leq i\leq c$. Thus $\A\notin\mt{BCK}_{c}$.

\end{proof}

From this we see that $\mt{BCK}_c$ is a subquasivariety of $\mt{BCK}$ with quasiequational basis given by (BCK1)-(BCK5) together with $[x_1, x_2, \ldots, x_c]=0$. On the other hand, 

\begin{theorem}
$\mt{BCK}_c$ is a variety if and only if $c=1$. 
\end{theorem}

\begin{proof}
The nilpotent BCK-algebras with nilpotence class 1 are precisely the commutative BCK-algebras; this is known to be a variety \cite{yutani77}. 

Suppose $c\geq 2$. In \cite{wronski83}, Wro\'{n}ski gave an example of a BCK-algebra $\U$ with a congruence $\theta$ such that $\U/\theta$ is not a BCK-algebra. We will show this algebra $\U$ is nilpotent of class 2. 

We first define the algebra. Let $N=\{0,1,2,\ldots\}$, $A=\{a_n\,\mid\, n\in N\}$, and $B=\{b_n\,\mid\, n\in N\}$. The carrier set is $U=N\cup A\cup B$, and this becomes a BCK-algebra $\U$ under the following operation:
\begin{align*}
n\bdot m &= \max\{0,n-m\}\\
n\bdot a_m &=n\bdot b_m=0\\
a_n\bdot m &= a_{n+m}\\ 
b_n\bdot m &= b_{n+m}\\
a_n\bdot a_m &=b_n\bdot b_m=\max\{0,m-n\}\\
a_n\bdot b_m &=a_n\bdot a_{m+1}\\
b_n\bdot a_m &= b_n\bdot b_{m+1}
\end{align*} for all $n,m\in N$. The Hasse diagram for $\U$ is shown in Figure \ref{fig:U}.
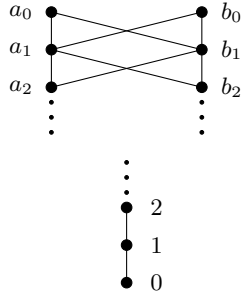
\begin{figure}[h]
\centering
\begin{tikzpicture}
\filldraw (0,0) circle (2pt);
\filldraw (0,.5) circle (2pt);
\filldraw (0,1) circle (2pt);
\filldraw (0,1.2) circle (.7pt);
\filldraw (0,1.4) circle (.7pt);
\filldraw (0,1.6) circle (.7pt);
\filldraw (-1,2) circle (.7pt);
\filldraw (-1,2.2) circle (.7pt);
\filldraw (-1,2.4) circle (.7pt);
\filldraw (1,2) circle (.7pt);
\filldraw (1,2.2) circle (.7pt);
\filldraw (1,2.4) circle (.7pt);

\filldraw(-1,2.6) circle (2pt);
\filldraw(-1,3.1) circle (2pt);
\filldraw(-1,3.6) circle (2pt);
\filldraw(1,2.6) circle (2pt);
\filldraw(1,3.1) circle (2pt);
\filldraw(1,3.6) circle (2pt);
\draw [-] (0,0) -- (0,1);
\draw [-] (-1,2.6) -- (-1,3.6);
\draw [-] (1,2.6) -- (1,3.6);
\draw [-] (-1,2.6) -- (1,3.1);
\draw [-] (1,2.6) -- (-1,3.1);
\draw [-] (-1,3.1) -- (1, 3.6);
\draw [-] (1,3.1) -- (-1, 3.6);
	\node at (.4,0) {\small 0};
	\node at (.4, .5) {\small 1};
	\node at (.4, 1) {\small 2};
	\node at (-1.4, 2.6) {\small $a_2$};
	\node at (-1.4, 3.1) {\small $a_1$};
	\node at (-1.4, 3.6) {\small $a_0$};
	\node at (1.4, 2.6) {\small $b_2$};
	\node at (1.4, 3.1) {\small $b_1$};
	\node at (1.4, 3.6) {\small $b_0$};
\end{tikzpicture}
\caption{Hasse diagram for $\U$}
\label{fig:U}
\end{figure} 

Some computation reveals that
\[[a_n, b_m]=[b_n,a_m]=\begin{cases}0&\text{ if $n<m$}\\1&\text{ if $n=m$}\\2&\text{ if $n>m$}\end{cases}\] and all other pseudocommutators are 0. The set $\{0,1,2\}$ is already a subalgebra and so $\U'=\{0,1,2\}$. Since all pseudocommutators involving 0, 1, or 2 are again 0, we have $\U_2=\{0\}$. Hence, $\U$ is nilpotent of class 2, and therefore $\U\in \mt{BCK}_c$ for any $c\geq 2$. But as mentioned, $\U$ has a homomorphic image which is not a BCK-algebra, which implies that $\mt{BCK}_c$ is not a variety for $c\geq 2$.

\end{proof}

In \cite{NS14}, it was observed that all BCK-algebras up to order five are solvable, and the authors conjectured that all finite BCK-algebras are solvable. Dudek resolved this conjecture in the affirmative in \cite{dudek16}. The same is true for nilpotency; in fact, we prove the more general result that any BCK-algebra of finite height is nilpotent. The proof strategy is similar to Dudek's.

\begin{lemma}[\hspace{1pt}\cite{dudek16}, Lemma 1(c)]\label{lem1}
If $x\neq 0$, then $[x,y]<x$ for any $y\in \A$. 
\end{lemma}

\begin{lemma}\label{A-M subalg}
Let $\A$ be a non-trivial BCK-algebra and let $M$ denote the subset of maximal elements. Then $\A\setminus M$ is a subalgebra of $\A$ and $\A'\leq \A\setminus M$. 
\end{lemma}

\begin{proof}
For the first claim, if $\A\setminus M=\{0\}$, we are done. So suppose $\A\setminus M\neq\{0\}$. Take $u,v\in \A\setminus M$. Then $u\bdot v\leq u$ which shows $u\bdot v$ is not a maximal element. Hence $u\bdot v\in \A\setminus M$. 

To show $\A'\leq \A\setminus M$, we will show that $M$ contains no pseudocommutators. To this end, first notice that $[0,y]=0\notin M$ for all $y\in \A$ since $\A$ is non-trivial. Then, for $x\neq 0$, we have $[x,y]<x$ for any $y\in\A$ by Lemma \ref{lem1} which implies $[x,y]$ cannot be a maximal element for any $y\in\A$. Thus, $[x,y]\notin M$ for any $x,y\in \A$, meaning $\A\setminus M$ contains all the pseudocommutators and therefore $\A'\leq \A\setminus M$.

\end{proof}

Define the {\bf height} of $\A$, denoted $\Ht(\A)$, to be the height of $\A$ as a poset; that is, $\Ht(\A)$ is the maximum cardinality of any chain in $\A$. Note that if $\A$ has finite height, then $\Ht(\A')\leq \Ht(\A\setminus M)=\Ht(\A)-1<\Ht(\A)$.

\begin{theorem}\label{fin = nil}
Any BCK-algebra $\A$ of finite height is nilpotent of class at most $\Ht(\A)$.
\end{theorem}

\begin{proof}
Suppose $\A$ is a BCK-algebra of finite height. Then $M\neq\emptyset$. If $\A'=\{0\}$ then $\A$ is nilpotent and we are done. 

Assume $\A'\neq\{0\}$. Let $M'\subseteq \A'$ denote the maximal elements of $\A'$. We claim that \[\bigl\{\,[x,y]\,\mid\,x\in\A', y\in\A\,\bigr\}\subseteq \A'\setminus M'\,.\] Since $\A$ has finite height and $\A'\neq\{0\}$, we have $M'\neq\emptyset$ and $0\notin M'$, and thus $[0,y]=0\notin M'$ for any $y\in\A$. Now take $x\neq 0$ in $\A'$. By Lemma \ref{lem1} we have $[x,y]<x$ and so $[x,y]$ cannot be maximal in $\A'$ for any $y\in\A$. This proves the claim which now implies that $\A_{(2)}=[\A',\A] \leq \A'\setminus M'$. Since $M'\neq\emptyset$, we see that $\Ht(\A_{(2)})<\Ht(\A')$. By an induction we have
\[\Ht(\A)>\Ht(\A')>\Ht(\A_{(2)})>\Ht(\A_{(3)})>\cdots\] but this means $\Ht(\A_{(k)})=1$ for some $k$ since $\A$ has finite height. Then $\A_{(k)}=\{0\}$, meaning $\A$ is nilpotent of class $k\leq\Ht(\A)$.

\end{proof}

\begin{corollary}
Any finite BCK-algebra is nilpotent.
\end{corollary}

%
%
%
%
%


\end{document}